\documentclass[a4paper, 10pt]{article}

\usepackage{exscale}
\usepackage[centertags]{amsmath}
\usepackage[T1]{fontenc}
\usepackage{amssymb}
\usepackage{amsthm}
\usepackage{varioref}
\usepackage[curve,matrix,arrow,ps]{xy}
\usepackage{epsfig}

\newtheorem{definition}{Definition}
\newtheorem{satz}{Proposition}
\newtheorem{theorem}{Theorem}
\newtheorem{proposition}{Proposition}
\newtheorem{lemma}{Lemma}
\newtheorem{korollar}{Corollary}

\newcommand{\nd}{\noindent}

\newcommand{\dI}{{\mathbf 1}}
\newcommand{\real}{{\mathbb R}}
\newcommand{\dC}{{\mathbb C}}

\newcommand{\dN}{{\mathbb N}}
\newcommand{\dK}{{\mathbb K}}

\newcommand{\dH}{{\mathbb H}}
\newcommand{\dL}{{\mathbb L}}

\newcommand{\dR}{{\mathbb R}}

\newcommand{\cC}{\mathcal{C}}
\newcommand{\cD}{\mathcal{D}}
\newcommand{\cE}{\mathcal{E}}
\newcommand{\cF}{\mathcal{F}}
\newcommand{\cG}{\mathcal{G}}
\newcommand{\cH}{\mathcal{H}}
\newcommand{\cI}{\mathcal{I}}

\newcommand{\cK}{\mathcal{K}}
\newcommand{\cL}{\mathcal{L}}
\newcommand{\cM}{\mathcal{M}}
\newcommand{\cN}{\mathcal{N}}
\newcommand{\cO}{\mathcal{O}}
\newcommand{\cP}{\mathcal{P}}

\newcommand{\cT}{\mathcal{T}}

\newcommand{\D}{\displaystyle}

\newcommand{\SC}{\scriptstyle}

\title{Deformation of singular lagrangian subvarieties}
\author{Duco van Straten and Christian Sevenheck}
\date{February 10, 2000}

\begin{document}
\maketitle
\begin{abstract}
We investigate deformations of lagrangian manifolds with
singularities. We introduce a complex similar to the
\emph{de Rham}-complex whose cohomology calculates deformation
spaces. Examples of singular lagrangian varieties are presented
and deformations are calculated explicitly.
\end{abstract}
\renewcommand{\thefootnote}{}
\section{Introduction}

In this paper, we develop some ideas of a deformation theory of
singular lagrangian subvarieties. Lagrangian submanifolds are
quite fundamental objects, so in a sense it is natural to extend
the study of them to a larger class of objects which are allowed
to have singularities. This has been done by Arnold, Givental and
others (\cite{Giv3}). However, not much is known on the behavior
of lagrangian singularities under deformations. The aim of this
article is to describe the spaces of infinitesimal deformations
and obstructions of a lagrangian subvariety and to perform
calculations for some concrete examples. It turns out that the
lagrangian property of a space has a strong influence on
its deformations, e.g., there are examples of spaces $X$ with
$dim(T^1_X) =  \infty$, which have nevertheless a versal
deformation space for the lagrangian deformations.

In the sequel, we will consider the following situation: Let $M$
be a $2n$-dimen\-sio\-nal symplectic manifold over $\dK = \real$
or $\dK=\dC$ (that is, a $C^\infty$ or complex analytic manifold
of real resp. complex dimension $2n$ endowed with a closed,
non-degenerated 2-form $\omega$, holomorphic in the second case)
and $L$ a reduced analytic subspace of dimension $n$, given by an
involutive ideal sheaf $\cI$, i.e. an ideal sheaf satisfying
$\{\cI,\cI\} \subset \cI$ where $\{\,,\}$ denotes the Poisson
bracket corresponding to $\omega$. This condition ensures that $L$
is a lagrangian submanifold in a neighborhood of each of its
smooth points. A lagrangian deformation of $L$ will be a
deformation in the usual sense (a flat family $L_S \rightarrow S$)
with the additional condition that all fibers are lagrangian
subvarieties of $M$. More precisely, we will call a diagram
$$
 \xymatrix@!0{
  & M \ar@(d,ur)[dddl] \ar@{^{(}->}[rr] && M \times S \ar@(d,ur)[dddl]   \\
  L \ar@{^{(}->}[ur]\ar[rr] \ar[dd]  &  & L_S \ar@{^{(}->}[ur]\ar[dd] \\
  \\
  \{*\} \ar[rr]  &  & S        }
$$
a lagrangian deformation of $L$ iff $L_S \rightarrow S$ is flat
and $\{\cI_S,\cI_S\}_S \subset \cI_S$. Here $\cI_S$ is the ideal
sheaf defining $L_S$ in $M\times S$ and $\{\;,\}_S$ is the Poisson
structure defined on $M \times S$ by the (degenerate) form
$\omega_S:=p^*\omega$, $p: M\times S \rightarrow M$ being the canonical
projection. This definition can be formalized using the
language of deformation functors (see \cite{Dipl} and \cite{Schlessinger}). This
more formal approach yields the definition of morphisms of
deformations, in particular, two deformations $L_S \subset M\times S$
and $L'_T \subset N\times T$ are called equivalent iff there is an
fibrewise isomorphism $F:M\times S\rightarrow N\times T$ satisfying $F^*\omega_T
=\omega_S$. Such an $F$ comes from a symplectic diffeomorphism
$f:M \rightarrow N$ and in case that $M$ is simply connected
(which we will suppose from now on), $f$ is induced by an
hamiltonian vector field, see lemma \ref{lemSympHamil}.

The tangent space to the functor of lagrangian deformations of
$L$ (that is, the space of lagrangian deformations of $L$
over $Spec(\dK[\epsilon])$ up to those induced by hamiltonian vector
fields of the ambient manifold) will be denoted by
$LT^1_L$. However, we will focus our attention to the local case mainly,
that is, we will study the \emph{sheaf} $\cL\!\cT^1_L$ of lagrangian deformations
of $L$. For lagrangian submanifolds, it follows from \cite{Voisin}, that the versal
deformation space is smooth, i.e., deformation of such objects are unobstructed.
This is probably not true in the singular case, although an example has not been
found yet. See theorem \ref{theMain} for further details.

\textbf{Acknowledgements}: We would like to thank A.~Givental who suggested
to investigate the deformation theory of lagrangian singularities in december
1992.

\section{The complex $\cC^\bullet$}

We start with a slightly more general situation: Let $\cI \subset \cO_M$
be an involutive ideal sheaf, $\cO_L$ the structure sheaf of the
subvariety $L$ described by $\cI$ and denote by $\cL:=\cI/\cI^2$ the
conormal sheaf. The formula $\{\cI^i,\cI^j\}\subset\cI^{i+j-1}$,
which can be easily verified, shows that there are well-defined
operations
$$
\begin{array}{rcl}
\cL \times \cO_L  &\longrightarrow &\cO_L \\
(g,f)&\longmapsto&\{g,f\}
\end{array}
\;\;\;
\mbox{and}
\;\;\;
\begin{array}{rcl}
\cL \times \cL &\longrightarrow &\cL \\
(g,h)&\longmapsto&\{g,h\}
\end{array}
$$
compatible in the sense that $\{g,f\cdot h\}=\{g,f\}h+f\{g,h\}$. This implies
that the first operation can be rewritten as a morphism
$$
\cL \rightarrow {\cD}\!er(\cO_L,\cO_L) = \Theta_L
$$
One says that $\cL$ is a \emph{Lie algebroid} (for details on Lie algebroids,
see \cite{Mackenzie}).
\begin{definition}
Let $\cC^p_L$ be the following $\cO_L$-module
$$
\cC^p_L := {\cH}om_{\cO_L}\left(\bigwedge^p \cL, \cO_L\right)
$$
and define a differential:
$$
\begin{array}{l}
\left(\delta\left(\phi\right)\right)\left(h_1 \wedge \ldots \wedge h_{p+1}\right) := \\
  \sum_{i=1}^{p+1} \left(-1\right)^i
  \left\{h_i,\phi\left(h_1\wedge\ldots \wedge \widehat{h}_i \wedge \ldots
  h_{p+1}\right)\right\}
\\
 +  \sum\limits_{1 \leq i < j \leq p+1}\left(-1\right)^{i+j-1}
\phi\left(\left\{h_i,h_j\right\} \wedge  h_1  \wedge \ldots \wedge \widehat{h}_i
\wedge \ldots \wedge \widehat{h}_j \wedge \ldots \wedge h_{p+1} \right)
\end{array}
$$
\end{definition}
It is a straightforward computation to check that
$\delta \circ \delta = 0$, so we get indeed a complex.
Following \cite{Mackenzie}, it is called the standard
complex for the Lie algebroid $\cL$. Remark that $\cC^0=\cO_L$
and $\cC^1={\cH}\!om_{\cO_L}(\cI/\cI^2,\cO_L)=:\cN_L$, the normal sheaf
of $\cI$ in $\cO_M$. For the definition of $\delta$, the fact that
$\cI$ is involutive is essential: the second term would not make
sense otherwise.

\nd
We may define a product on the complex $(\cC^\bullet,\delta)$:
\begin{eqnarray*}
\cC^p \times \cC^q & \longrightarrow & \cC^{p+q} \\
(\Phi,\Psi) & \longmapsto & \Phi \wedge \Psi
\end{eqnarray*}
with
$$
\begin{array}{lll}
(\Phi \wedge \Psi)(f_1 \wedge \ldots \wedge f_{p+q})
 =  & & \\
\\
\D
\sum\limits_{
\begin{array}{c}
\SC
  I \coprod J = \{1,\ldots,n\}\\
\SC
  i_1 < \ldots < i_p \\
\SC
  j_1 < \ldots < j_q
\end{array}
} sgn(I,J) \cdot
\Phi(f_{i_1}\wedge\ldots\wedge f_{i_p}) \cdot \Psi(f_{j_1}\wedge\ldots\wedge f_{j_q})
&&
\end{array}
$$
The sign is defined as
$$
sgn(I,J) := sgn
\binom
{1,\dotfill,p+q}
{i_1,\ldots,i_p,j_1,\ldots,j_q}
$$
\begin{satz}
\label{propDGA}
Let $\Phi \in \cC^p$, $\Psi \in \cC^q$ et $\Gamma \in \cC^r$.
Then we have
\begin{enumerate}
\item
$
\Phi \wedge \Psi = (-1)^{deg(\Phi) \cdot deg(\Psi)} \cdot \Psi \wedge \Phi
$
\item
$
(\Phi \wedge \Psi) \wedge \Gamma = \Phi \wedge (\Psi \wedge\Gamma)
$
\item
$
\delta(\Phi \wedge \Psi) =
\delta(\Phi) \wedge \Psi  + (-1)^{deg(\Phi)} \cdot \Phi \wedge
\delta(\Psi)
$
\end{enumerate}
\end{satz}
\begin{proof}
The first two points are trivial, while the third has to be
checked by an explicit calculation.
\end{proof}

Note that the last proposition says that $(\cC_L^\bullet,\delta,\wedge)$ is
a \emph{differential graded algebra}, furthermore, we have $\cC^0_L = \cO_L
=\Omega^0_L$. As one might hope, there is indeed a tight
connection between $\Omega^\bullet_L$ and $\cC^\bullet_L$.
\begin{satz}
Suppose that $L$ is lagrangian. Then there exists a morphism
$J:\Omega^1_L \rightarrow \cC^1_L$ which is an isomorphism
outside the singular locus of $L$.
\end{satz}
\begin{proof}
On a symplectic manifold, there is a canonical isomorphism $\beta$ between
vector fields and one forms, given by $\beta(V) := i_V \omega$.
On the other hand, for each analytic subspace $L \subset M$ we
have two exact sequences, dual to each other, namely, the conormal
and the normal sequence, thus, there is the following diagram:
$$
\xymatrix@C=0.5cm{
&{\cL} \ar[r] & \Omega^1_M \otimes \cO_L \ar[r] \ar[d]^{\alpha:=\beta^{-1}} & \Omega^1_L \ar[r] & 0 \\
0 \ar[r] &\Theta_L \ar[r] & \Theta_M \otimes \cO_L \ar[r] & {\cN_L} \ar[r] & {\cT^1_L} \ar[r] &0 }
$$
Now the fundamental fact is that this diagram can be completed:
the morphism $\cL \rightarrow \Theta_L$ above commutes with
$\alpha$, so we have
\begin{equation}
\label{eqAlpha'}
\xymatrix@C=0.5cm{
{\cL} \ar[r] \ar[d]_{\alpha'} & \Omega^1_M \otimes \cO_L \ar[d]^{\alpha} \\
\Theta_L \ar[r] & \Theta_M \otimes \cO_L }
\end{equation}
Note that the image of an element $g\in \cL$ under $\alpha'$ is just the hamiltonian
vector field $H_g$. The morphisms $J:\Omega^1_L
\rightarrow \cC^1_L=\cN_L$ we are looking for can now be defined as
the map induced by $\alpha$, explicitly
$$
J(df) = \left(g \mapsto \{f,g\}\right)
$$
To see that $J$ is an isomorphism near a smooth point of $L$ it
will be sufficient to prove this for the map $\alpha'$ (because
at smooth points $x$ we have $\cT^1_{(L,x)} = 0$ and the map
$\cL_x \rightarrow \Omega^1_{(L,x)} \otimes \cO_{L,x}$ is
injective). So assume the sheaves $\cL$, $\Omega^1_L$, and
$\Theta_L$ to be defined in a neighborhood of a smooth point
which means that they all become locally free.
$\cL$ then has to be identified with the conormal bundle.
To prove that $\alpha'$ is an isomorphism, we will construct an
inverse. First note that, by the fact that $L$ is
coisotropic, the morphism $\beta: {\Theta_M}_{|L} \rightarrow
{\Omega^1_M}_{|L}$ actually sends an element of $\Theta_L$
to a form vanishing on all vectors tangent to $L$. So the
restriction of $\beta$ to $\Theta_L$ defines a morphism $\beta':
\Theta_L \rightarrow \cL$. The situation is as follows:
$$
\xymatrix@C=0.5cm{
0 \ar[r] &\Theta_L \ar[r] \ar[d]^{\beta'} & \Theta_M \otimes \cO_L \ar[r]
\ar[d]^{\beta} & {\cN_L} \ar[r] &0 \\
0 \ar[r] &{\cL} \ar[r] & \Omega^1_M \otimes \cO_L \ar[r] & \Omega^1_L \ar[r] & 0}
$$
A diagram chase shows that $\beta'$ is injective. On the other
hand, we have $dim(\cL) = dim(\Theta_L)$, as $L$ is lagrangian. So
$\beta'$ is an isomorphism and the inverse of $\alpha'$.
\end{proof}
\begin{korollar}
\label{corJDach}
The morphism $J:\Omega^1_L\rightarrow\cC^1_L$ can be extended to
a morphism of DGA's
$$
J: (\Omega_L^\bullet,d,\wedge) \longrightarrow  (\cC_L^\bullet,\delta,\wedge)
$$
which is an isomorphism at smooth points of $L$.
\end{korollar}
\begin{proof}
Set
$$
J(\omega_1\wedge\ldots\wedge\omega_p):=
J(\omega_1)\wedge\ldots\wedge J(\omega_p)
$$
where $\omega_i \in \Omega^1_L$. Then it is immediate that
$J$ is an isomorphism on $L_{reg}$. To prove that
$J \circ d = \delta \circ J$, it suffices to
check this in the lowest degrees, that is, we have to show that
the diagram
$$
\xymatrix{
  \Omega^0_L \ar[d] \ar[r]^{d} & \Omega^1_L \ar[d]^{J} \\
  {\cC^0_L} \ar[r]^{\delta} & {\cC^1_L}}
$$
commutes. This follows directly from $\Omega^0_L=\cC^0_L=\cO_L$.
\end{proof}
In the last section, we use the following elementary fact.
\begin{lemma}
\label{lemTors}
The kernel of $J$ is the complex ${\cT}\!ors(\Omega_L^\bullet)$
consisting of the torsion subsheaves of $\Omega^p_L$.
\end{lemma}
\begin{proof}
We have ${\cT}\!ors(\Omega^\bullet_L) \subset {\cK}\!er(J)$
as $\cC^\bullet_L$ is torsion free. On the other hand, the kernel
is supported on the singular locus of $L$, so it must be a torsion
sheaf, hence ${\cK}\!er(J) \subset
{\cT}\!ors(\Omega^\bullet_L)$.
\end{proof}
\paragraph{Remark:} Although the definition of the modules
$\cC_L^p$ involves the ideal $\cI$, they are probably intrinsic.
This is at least clear in some special cases as the following
lemma shows.
\begin{lemma}
Suppose $L$ to be Cohen-Macaulay and regular in codimension one.
Then there is an isomorphism
$$
(\Omega_L^p)^{**} \stackrel{\cong}{\longrightarrow} \cC_L^p
$$
where for an $\cO_L$-module $\cF$, $\cF^*$ denotes ${\cH}\!om_{\cO_L}
(\cF,\cO_L)$.
\end{lemma}
\begin{proof}
We will make use of the following fact: Let $\cF$ be an $\cO_L$-module
of type $\cG^*$, then $\cF$ is reflexive, i.e. $\cF^{**}=\cF$. The morphism
$h: (\Omega_L^p)^{**}\rightarrow \cC_L^p$ we are looking is obtained by
dualizing twice the morphism $J:\Omega_L^p \rightarrow \cC_L^p$, this yields
$J^{**}:(\Omega_L^p)^{**} \rightarrow (\cC_L^p)^{**} = \cC_L^p$ as $\cC^p_L$ is of
type ${\cH}\!om(-,\cO_L)$. Clearly, $h$ is an isomorphism on the regular locus.
We have an exact sequence
$$
0 \longrightarrow \cK
  \longrightarrow (\Omega_L^p)^{**}
  \stackrel{h}{\longrightarrow} \cC_L^p
  \longrightarrow \cG
  \longrightarrow 0
$$
where $\cK$ and $\cG$ are the kernel resp. cokernel sheaves of the map
$h$. This sequence can be split
$$
\begin{array}{c}
0 \longrightarrow \cK
  \longrightarrow (\Omega_L^p)^{**}
  \longrightarrow \cH  \longrightarrow 0\\
0 \longrightarrow \cH
  \longrightarrow \cC_L^p
  \longrightarrow \cG \longrightarrow 0
\end{array}
$$
with $\cH = {\cI}\!m(h)$. Applying ${\cH}\!om_{\cO_L}(-,\cO_L)$ yields
$$
\begin{array}{c}
0 \longrightarrow \cH^*
  \longrightarrow ((\Omega_L^p)^{**})^*
  \longrightarrow \cK^*  \\
0 \longrightarrow \cG^*
  \longrightarrow (\cC_L^p)^*
  \longrightarrow \cH^* \longrightarrow
  {\cE}\!xt^1(\cG, \cO_L)
\end{array}
$$
Now we use the lemma of Ischebeck (see \cite{Matsumura}): Given a
local ring $R$, two $R$-modules $M$ and $N$ with $k=dim(M)$ and
$r=depth(N)$, then for all $p<r-k$, the modules $Ext^p(M,N)$
vanish. It follows that $\cK^* = \cG^* = {\cE}\!xt^1(\cG, \cO_L) =
0$, so we have $((\Omega_L^p)^{**})^*=(\cC_L^p)^*$. Then obviously
$((\Omega_L^1)^{**})^{**}=(\cC_L^1)^{**}$ and by the argument above
$(\Omega_L^1)^{**}=\cC_L^1$ so the map $h$ is an isomorphism.
\end{proof}

\section{Deformations}

Recall that the space of infinitesimal embedded deformations of an analytic
algebra $R$, given as $R=S/I$ where $S$ is the ring of convergent
power series, is equal to the normal module of $I$ in $S$, i.e.
$Hom_R(I/I^2,R)$. Dividing out trivial deformations gives the
space $T^1_R$, defined by the sequence
$$
\xymatrix@C=0.5cm{
0 \ar[r] & Hom_R(\Omega^1_{R},R) \ar[r] & Hom_S(\Omega^1_S,S)\, \widehat{\otimes} \, R \ar[r] &
Hom_R(I/I^2,R) \ar[r] & T^1_R \ar[r] & 0 }
$$
On the other hand, the deformations of a
manifold $X$ over $Spec(\dK[\epsilon]/(\epsilon^2))$ are
parameterized by $H^1(X,\Theta_X)$. The cotangent complex is a tool
to handle these two special cases in an integrated manner:
infinitesimal deformations of an analytic space $L$ are in
bijection with $\dH^1(\dL_X)$. It seems that the complex $\cC_L^\bullet$
has to be seen as a first approximation to an equivalent
for the cotangent complex in the lagrangian
context. More precisely, the following holds:
\begin{theorem}
\label{theMain}
The first three cohomology sheaves of $\cC^\bullet_L$ are
\begin{itemize}
\item
$\cH^0(\cC_L^\bullet)=\dK_L$.
\item
$\cH^1(\cC_L^\bullet)=\cL\!\cT^1_L$.
\item
$\cH^2(\cC_L^\bullet)=\cL\!\cT^2_L$. This symbol denotes the
\textbf{lagrangian obstructions}, that is, $\cL\!\cT^2_L$ is
the sheaf of obstructions to extend a lagrangian deformation to
higher order regardless whether it can be extended as a flat
deformation.
\end{itemize}
\end{theorem}
\nd
The proof of the following preliminary lemma can be found in \cite{Banyaga}.
\begin{lemma}
\label{lemSympHamil}
If $H^1(M, \dK)=0$, then each diffeomorphism $f:M \rightarrow M$ satisfying $f^*\omega =
\omega$ is the time $1$ map of a flow $\varphi_t$ of a hamiltonian vector
field $H_g$ for some function $g$ on M.
\end{lemma}
\begin{proof}[Proof of the theorem]
$\cH^0(\cC^\bullet_L)$ equals ${\cK}\!er(\delta:\cO_L\rightarrow\cC^1_L)$. Take an element
$f$ of ${\cK}\!er(\delta)$. Then $\{f,g\} \in \cI$ for all
$g \in \cI$. If $f$ is not a constant, then the ideal $(\cI, f)$ is
strictly larger than $\cI$, not the whole ring and still
involutive. This is a contradiction to the fact that $L$ is lagrangian,
which means that $\cI$ is maximal under all involutive ideals. So
the kernel must be the constant sheaf.

To prove that $\cH^1(\cC_L^\bullet)=\cL\!\cT^1_L$, two things have to be
checked: As $\cC_L^1=\cN_L$, we must first identify the elements of
${\cK}\!er(\delta^1:\cC^1_L\rightarrow\cC^2_L)$ with the flat
\textbf{lagrangian} deformations. Then we have to show that the
image of $\delta^0:\cO_L\rightarrow\cC^1_L$ are the trivial
deformations. But this is easy, because for $f\in\cO_L$,
$\delta(f)$ acts as $H_f$, thus inducing a trivial deformation.
Furthermore, by lemma \ref{lemSympHamil}, of all deformations coming from vector fields on $M$,
only those induced by hamiltonian vector fields are trivial in the
lagrangian sense. Now we choose an open set $U \subset L$ and
sections $(f_1,\ldots,f_k)$ generating $\cI(U)$. Take an element
$\Phi \in {\cK}\!er(\delta^1)$, which means that
$$
\phi\left(\{g,h\}\right)-\left\{g,\phi(h)\right\}-
\left\{\phi(g),h\right\} = 0
$$
for all $f,g \in \cI/\cI^2$. Then $\Phi$ corresponds to the
deformation given by
$$
\widetilde{\cI}=
\left(f_1+\epsilon \phi(f_1),\ldots,
f_k+\epsilon \phi(f_k)\right)
$$
The ideal $\widetilde{\cI}$ is involutive iff for any two elements
$f+\epsilon \phi(f), g+\epsilon \phi(g)$, we have $\left\{f+\epsilon
\phi(f), g+\epsilon \phi(g)\right\} \in \widetilde{\cI}$, which is
equivalent to
$$
F:=\left\{f,g\right\} + \epsilon \left( \{ f,\phi (g) \}
+\{ \phi (f),g\}\right) \in \widetilde{\cI}
$$
Consider $G:=\left\{f,g\right\} + \epsilon
\phi\left(\{f,g\}\right)$, which is an element of
$\widetilde{\cI}$, so the condition $F \in \widetilde{\cI}$ is
equivalent to $F-G \in \widetilde{\cI}$, that is
$$
\left\{ f,\phi \left(g\right) \right\}
+\left\{ \phi \left(f\right),g\right\}
-\phi\left(\left\{f,g\right\}\right) \in \cI
$$
This means exactly that $\phi \in {\cK}\!er(\delta^1)$.

In order to interpret the second cohomology group, we define the
bilinear mapping
$$
\begin{array}{rcl}
\widetilde{ob} : \cC^1_L \times \cC^1_L & \longrightarrow & \cC^2_L \\
(\Phi,\Psi) & \longmapsto & \left(g \wedge h \mapsto \{\Phi(g),\Psi(h)\}\right)
\end{array}
$$
In this way we get a quadratic form $ob(\Phi):=\widetilde{ob}(\Phi,\Phi)$.
It can be immediately verified that this induces a map $ob:\cH^1(\cC^\bullet_L)
\rightarrow \cH^2(\cC^\bullet_L)$. We will now prove the following:
Given a lagrangian deformation $\Phi \in \cL\!\cT^1_L$. Then there is a
lift to second order defining an involutive ideal iff $ob(\Phi)=0 \in \cL\!\cT^2_L$.
The last condition is equivalent to the existence of $\Psi \in
\cL\!\cT^1_L$ with $ob(\Phi) = \delta(\Psi)$, i.e.
$$
\left\{\Phi(f),\Phi(g)\right\}
=
 \Psi\left(\{f,g\}\right)
-\left\{f,\Psi(g)\right\}
-\left\{\Psi(f),g\right\}
\;\;\;\;\;\;\forall f,g \in \cL
$$
But this means that the following ideal is involutive.
$$
J=(f_1+\epsilon\Phi(f_1)+\epsilon^2\Psi(f_1), \ldots,
f_k+\epsilon\Phi(f_k)+\epsilon^2\Psi(f_k))
$$
\end{proof}
\paragraph{Remark:} The fact that $\cL\!\cT^2_L$ is not the real
obstruction space make precise what was meant by saying that
complex $\cC^\bullet_L$ is a first approximation of the object we
are looking for: Hopefully, there is a modified version of this complex
whose cohomology gives, in complete analogy with the cotangent
complex, the spaces $T^1$ and $T^2$ for \textbf{flat lagrangian}
deformations. On the other hand, it is perhaps not even necessary
to impose flatness as the involutivity condition implies that the
dimension cannot drop, see also \cite{Matsushita}.
\begin{korollar}
There is an exact sequence
$$
\xymatrix@C=0.5cm{
  0 \ar[r] & H^1(L,\dK_L) \ar[r] & \dH^1(\cC^\bullet_L) \ar[r] & H^0(L, \cL\!\cT^1_L)
  \ar[r] & H^2(L,\dK_L) \ar[r] & \dH^1(\cC^\bullet_L) }
$$
Furthermore, there are two special cases:
\begin{itemize}
\item
Let $L$ be a contractible space. Then $\dH^1(\cC^\bullet_L)=
H^0(L,\cL\!\cT^1_L)$ and in fact: $LT^1_L=H^0(L,\cL\!\cT^1_L)$.
\item
Let $L$ be Stein and smooth. Then it follows that $\dH^1(\cC^\bullet_L)=H^1(L,\dK_L)$
and the space of global deformations is indeed $LT^1_L=H^1(L,\dK_L).$
\end{itemize}
\end{korollar}
\begin{proof}
The first fact is just the definition of the sheaf $\cL\!\cT^1_L$.
In the second case, note that the space of embedded flat deformations is $H^0(L, \cN_L)$, where
$\cN_L$ is the normal bundle of $L$ in $M$. As $L$ is smooth, this happens to
be $H^0(L, \Omega^1_L)$, so each infinitesimal flat deformation corresponds to
globally defined one-form on $L$. It is closed iff the deformation is lagrangian
and the subspace of exact one-forms are deformations induced by hamiltonian vector fields
(isodrastic deformations, see \cite{Weinstein}), these are the trivial ones.
$L$ is assumed to be a Stein manifold, in this case the first \emph{de
Rham}-cohomology group is exactly $H^1(L, \dK_L)$.
\end{proof}
By analogy with the cotangent complex, the following
generalization is probably true although we did not check the
details.
\begin{proposition}
The space of infinitesimal lagrangian deformations of a complex
space $L$ which is a lagrangian subvariety of a symplectic
manifold $(M,\omega)$ is given by
$$
LT^1_L = \dH^1(\cC_L^\bullet)
$$
\end{proposition}

\section{Finiteness of the cohomology}

This section is devoted to the proof of the following result.
\begin{theorem}
\label{theFiniteness}
Let $L \subset M$ be a lagrangian subvariety as above. Assume that
the following condition is satisfied: Denote by $edim(p)$ the
embedding dimension of a point $p \in L$, that is $edim(p) :=
dim_\dK(\mathbf{m}_p/\mathbf{m}^2_p)$, where $\mathbf{m}_p$ is the
maximal ideal in the local ring $\cO_{(L,p)}$. Let $S^L_k$ be the
following set
$$
S^L_k := \{p\in L\,|\,edim(p)=2n-k\} \subset L
$$
for all $k \in \{0,\ldots,n\}$, then suppose that we have
\begin{equation*}
\label{condition_P}
dim(S^L_k) \leq k
\end{equation*}
for all $k$. Under this condition (which will be called ``condition P''),
all $\cH^i(\cC^\bullet_L)$ are constructible sheaves of
$\dK$-vector spaces with respect to the stratification given
by the $S^L_k$.
\end{theorem}
Before going into the details of the proof, we would like to
explain the meaning of the condition (\ref{condition_P}).
\begin{lemma}\label{lemZerlegung}
Let $p \in S^L_k \subset L$ with $k > 0$. Then the germ $(L,p)$ can be
decomposed into a product
$$
(L,p) = (L',p') \times (\dK,0)
$$
where $(L',p')$ is a germ of a lagrangian variety in the
symplectic space $\dK^{2n-2}$. Furthermore, we have
$p' \in S^{L'}_{k-1}$.
\end{lemma}
\begin{proof}
Let $x_1,\ldots,x_{2n}$ be coordinates of $M$ centered at $p$.
Then the fact that $edim(p)<2n$ implies that there are
coefficients $\alpha_i \in \cO_{L,p}$ such that the following
equation holds in $\cO_{L,p}$
$$
\sum_{i=1}^{2n} \alpha_i x_i + h = 0
$$
where $h$ is an element of $\cO_{L,p}$ vanishing at second order.
So we have an element in the ideal describing $(L,p)$ whose
derivative do not vanish. Then $(L,p)$ is fibred by the
hamiltonian flow of this function. Explicitly, we can make
an analytic change of coordinates, such that $\alpha_1=1$,
$\alpha_i=0$ for all $i>1$ and $h=0$. Than the ideal of $(L,p)$
is of the form $(x_1,f_1,\ldots,f_m)$ for some functions $f_i$
which are independent of the variable $x_{n+1}$ (provided that we have
chosen the symplectic form to be $\sum_{i=1}^n d x_i\wedge
d x_{n+i}$).
\end{proof}
According to the lemma, the set of points of the variety $L$ can
be divided into two classes, those with maximal embedding dimension
(these are the ``bad points'') and those (with $edim(p) < 2n$) at
which $L$ is decomposable. Condition P
implies that the bad points are isolated. As usual, the proof of
the theorem consists of two parts: First, we will show that the
cohomology sheaves are locally constant on the strata $S^L_k$.
This is an immediate consequence of the following lemma. Then it
suffices to show that all stalks of $\cH^p(\cC_L^\bullet)$ are finite-dimensional.
\begin{lemma}[Propagation of Deformations]
\label{lemDecomp}
Let
$$
(L,0) \subset (\dK^{2n},0)
$$
be a germ of a lagrangian subvariety which can be decomposed,
i.e., there is a germ $(L',0)$ (which is lagrangian in
$(\dK^{2n-2},0)$) such that $(L,0)=(L',0)\times(\dK,0)$. Denote by
$\pi:L \rightarrow L'$ the projection. Then there is a
quasi-isomorphism of sheaf complexes
$$
j: \pi^{-1}\cC^\bullet_{L'} \rightarrow \cC^\bullet_L
$$
\end{lemma}
\begin{proof}
The proof of lemma \ref{lemZerlegung} shows that the ideals
$I$ and $I'$ describing the two germs differ by exactly one
element whose differential do not vanish at the origin. This
implies that the conormal sheaves $\cL$ of $L$ and $\cL'$ of $L'$
are related by the formula $\cL = \pi^*\cL' \oplus \cO_L$. It
follows that
$$
\cC^p_L =
{\cH}\!om_{\cO_L}\left(\pi^*\bigwedge^p \cL',\cO_L\right)
\oplus
{\cH}\!om_{\cO_L}\left(\pi^*\bigwedge^{p-1} \cL',\cO_L\right)
$$
Now we have to describe the differential on $\cC^\bullet_L$. We
choose local Darboux coordinates $(p_1,\ldots,p_n,q_1,\ldots,q_n)$
on $\dK^{2n}$ and $(p_2,\ldots,p_n,q_2,\ldots,q_n)$ on
$\dK^{2n-2}$. Suppose that the two ideals are $I=(f_1,\ldots,f_m,p_1)$
and $I'=(f_1,\ldots,f_m,p_1,q_1)$ (if we consider $L'$ as embedded in
$\dK^{2n}$). Let $\Phi$ be an element of
$$
{\cH}\!om_{\cO_L} \left(\pi^*\bigwedge^p \cL',\cO_L\right)
$$
Then it can be written as a power series in $q_1$ with coefficients
in $\cC^\bullet_{L'}$. A direct calculation shows that the
differential on $\cC^\bullet_L$ is
$$
\begin{array}{ccccc}
\delta: &   \cC^p_L   & \longrightarrow & \cC^{p+1}_L \\
        & \sum\limits_{i=0}^\infty\left(\Phi,\Psi\right)q_1^i &      \mapsto    &
\sum\limits_{i=0}^\infty \left(\delta\Phi_i,\delta\Psi_i
+(-1)^{p+1}(i+1)\Phi_{i+1}\right)q_1^i
\end{array}
$$
It is clear that the morphism $j$ must be the obvious
inclusion
$$
{\cH}\!om_{\cO_{L'}}\left(\bigwedge^p \cL',\cO_{L'}\right)
\hookrightarrow
{\cH}\!om_{\cO_L}\left(\pi^*\bigwedge^p \cL',\cO_L\right)
\oplus
{\cH}\!om_{\cO_L}\left(\pi^*\bigwedge^{p-1} \cL',\cO_L\right)
$$
We will now show that the cokernel of this inclusion is acyclic.
Then it follows immediately that $j$ induces an isomorphism on
the cohomology. So let $\Gamma$ be an element of
${\cC}\!oker(j) \cap {\cK}\!er(\delta)$, that is,
$$
\Gamma=\sum_{i=1}^\infty(\Phi_i,\Psi_i)q_1^i+(0,\Psi_0)
$$
where $\delta\Phi_i=0$ and $\delta\Psi_i=(-1)^p(i+1)\Phi_{i+1}$ for
all $i$. But then $\Gamma$ vanishes in the cohomology because
it can be written as $\Gamma = \delta\Lambda$ with
$$
\Lambda :=
\sum_{i=1}^\infty\left(\frac{(-1)^p\Psi_{i-1}}{i},0\right)q_1^i
\in \cC_L^{p-1}
$$
\end{proof}
\begin{korollar}
\label{corTransConst}
We have isomorphisms of sheaves
$$
\pi^{-1}\cH^i(\cC^\bullet_{L'}) \cong \cH^i(\cC^\bullet_L)
$$
\end{korollar}
\begin{proof}
This is obvious since $\pi^{-1}$ is an exact functor.
\end{proof}
Let $p \in S^L_k$ be a point at which $L$ is decomposable, i.e.
$k>0$. By induction, we find a neighborhood $U \subset L$ of $p$
such there is an analytic isomorphism $h:U \stackrel{\cong}{\longrightarrow}
Z \times B_\epsilon(0)^k$, where $Z$ is lagrangian in
$\dK^{2(n-k)}$, $B_\epsilon(0):=\{z\in\dK\,|\,|z|<\epsilon\}$ and
each $q \in U \cap S^L_l$ corresponds via $h$ to a point $(q',b)
\in Z \times B(\epsilon)^k$ with $q'\in S^Z_{l-k}$. In particular,
the image of $U \cap S^L_k$ under $h$ is $(\{pt\},B(\epsilon)^k)$,
so by the last corollary, $\cH^p(\cC_L^\bullet)$ is
constant on $U\cap S^L_k$.

It remains to show that the stalks of the cohomology are
finite-dimensional. Again by corollary \ref{corTransConst}, this
is done once we have shown it for points with maximal embedding
dimension. We will use a method developed in
\cite{Buchweitz}. In this paper, the following situation is
considered. Let $f:X \rightarrow S$ be a morphism of complex
spaces (with $dim(S)=1$) and $\cK^\bullet$ a certain sheaf complex on $X$. Then,
under suitable conditions, the relative hypercohomology $\dR^i\!
f_* \cK^\bullet$ are coherent sheaves of $\cO_S$-modules. The
proof of this theorems relies on a functional analytic
argument of Kiehl and Verdier (see \cite{Kiehl} or \cite{Douady})
which states, roughly speaking, that if the mapping induced on the
complex of sections of $\cK^\bullet$ by a small shrinking of the
open set (over which the sections are taken) is a quasi-isomorphism,
then the hypercohomology groups are finite dimensional vector
spaces. We are going to use this result in the form of \cite{Straten}.
\begin{lemma}
\label{lemBadPoints}
Let $(L,p) \in (\dK^{2n},0)$ be a germ of a lagrangian variety
satisfying condition P which is indecomposable
at $p$. Then the stalk $\cH^i(\cC_L)_p$ is a finite-dimensional
$\dK$-vector space.
\end{lemma}
\begin{proof}
Choose a representative $V$ for the the germ such that
$edim(q) = 2n$ iff $q=p$ for all points $q \in V$. We refer the
reader to theorem 1 in \cite{Straten}. We do not
consider a relative situation here, so the map $f:X\rightarrow S$
in this theorem is replaced by $V\rightarrow\{0\}$ (Obviously, $V$
can be chosen such that this map is a standard representative of
the germ $(L,p)$ in the sense of definition 1
in \cite{Straten}, i.e., $V=L\cap B_\epsilon(p)$). The complex of sheaves in the theorem is the
complex $\cC^\bullet_L$, which satisfies the first two
properties ($C^p_L$ is $\cO_L$-coherent and the differential is
$\dK$-linear). Our task is to verify the third axiom, that is, we
have to find a vector field of class $C^\infty$ such that
$\cH^p(\cC_L^\bullet)$ is \emph{transversally constant} (see
definition 2 in \cite{Straten}), this will be done in corollary
\vref{lemGlobVect}. Now the proof of the theorem shows
that there is a smaller neighborhood $V_1$ of $p$ such that
$\Gamma(V, \cH^p(\cC^p_L)) = \Gamma(V_1, \cH^p(\cC^p_L))$. This gives
the result by using \cite{Kiehl} in the same way as in
\cite{Straten} or \cite{Buchweitz}.
\end{proof}
\begin{lemma}
Let $q \in V \cap S^L_k$ with $k>0$. Then there is a $C^\infty$-vectorfield
in a neighborhood $W$ of $q$ in $M$, tangent to $V \cap S^L_k$ and transversal to
$\partial B_\epsilon(p)$.
\end{lemma}
\begin{proof}
It follows from lemma \ref{lemZerlegung} that there exist $k$
linear independent hamiltonian vector fields on $M$ which respects
the stratum $S^k_L$. Now we have to distinguish the cases
$\dK=\real$ and $\dK=\dC$, in the first one, since $S^k_L$ is of
real dimension $k$ and since the intersection of $L$ and
$B_\epsilon(p)$ was transversal, it follows immediately that we
can find a linear combination of theses $C^\infty$-fields which is
transversal to $B_\epsilon(p)$. The same is true in the
complex case, here we have $k$ independent hamiltonian fields $\eta_1,
\ldots, \eta_k$ which are \emph{holomorphic}. As the holomorphic tangent space
at each point is canonically isomorphic (over $\real$) to the real one,
we get $2k$ linear independent $C^\infty$-fields by applying this isomorphism
to $\eta_1,\ldots,\eta_k,i \eta_1,\ldots,i \eta_k$. These can be used to find a
field transversal to $\partial B_\epsilon(p)$.
\end{proof}
\begin{korollar}
\label{lemGlobVect}
There is a $C^\infty$-vectorfield $\vartheta$ on a neighborhood
$U$ of $\partial B_\epsilon(p)$ in $M$ such that
$\cH^p(\cC^\infty_L)$ is transversally constant with
respect to $U$ and $\vartheta$.
\end{korollar}
\begin{proof}
Set $U:=(V\,\backslash\{p\})^\circ$. Then the last lemma
yields a covering $U_i$ of $U$ and vector fields $\vartheta_i$
defined in a neighborhood of $U_i$ in $M$. Chose a partition of unity
subordinate to this covering to obtain a field on $U$ which is
still transversal to $\partial B_\epsilon$. For each point $q\in
U$, which is contained in some stratum $S^L_k$, $\vartheta$ is
necessarily tangent to $S^L_k$, so the cohomology sheaves are
constant on the local integral curves of $\vartheta$.
\end{proof}
\paragraph{\textbf{Remark:}} By the \emph{Riemann-Hilbert-correspondence}
(see \cite{Björk}), the complex $\cH^\bullet:=\cH(\cC_L^\bullet)$, viewed as an object
of $\cD^b_c(\dK_M)$ (the derived category of constructible sheaves of $\dK$-vector spaces
on $M$) corresponds via the \emph{de Rham}-functor to a unique complex of
coherent $\cD_M$-modules with regular holonomic cohomology supported on $L$
(i.e., an object of $\cD^b_{\textup{r.h.}}(\mu_L(\cD_M))$).
\begin{lemma}
The complex $\cH^\bullet$ satisfies the first perversity condition, that is, the following
inequality holds.
$$
\textup{dim}(\cH^i(\cC^\bullet_L)) \leq n-i
$$
\end{lemma}
\begin{proof}
Let $p \in S^L_k$. Then $(L,p) = (L',p')\times (\dK^k,0)$ and
$\cH^i(\cC^\bullet_L)_p=\cH^i(\cC^\bullet_{L'})_{p'}$. But $dim(L')
\leq n - k$, so $\cH^i(\cC^\bullet_{L'})_{p'} = 0$ for all $i > n-k$.
\end{proof}
In case that the second perversity condition is also satisfied, the $\cH^i$'s are the
$\emph{de Rham}$-cohomology modules of some $\cD_M$-module supported on $L$. The following
consideration gives more evidence that the complex $\cC_L^\bullet$ is closely
related to $\cD$-module theory: Every complex manifold is lagrangian in its own cotangent bundle.
Consider \emph{Spencer's} complex, which is a resolution of $\cO_X$ as a $\cD_X$-module, explicitly:
$$
Sp(\cO_X)^\bullet:\hspace*{0.5cm}\ldots
\rightarrow \cD_X \otimes_{\cO_X} \Theta^{p+1}_X
\rightarrow \cD_X \otimes_{\cO_X} \Theta^p_X
\rightarrow \ldots \cD_X \rightarrow \cO_X \rightarrow 0
$$
The \emph{de Rham}-complex of $\cD_X$-module $\cM$ is obtained as
$$
DR(M) := {\cH}\!om_{\cD_X}(Sp(\cO_X)^\bullet, \cM)
$$
If we define a generalized version of the complex $\cC_L^\bullet$ as
$$
\cC^p_L(\cM) := {\cH}\!om_{\cO_L}\left(\bigwedge^p \cL, \cM\right)
$$
for some module $\cM$ over the Lie algebroid $\cL$, then $\cC^p_X(\cM)$ (for $X$ lagrangian
in $T^*\!X$) is exactly the \emph{de Rham}-complex of the $\cD_X$-Module
$\cM$.

\section{Examples and results}

In this section we will describe some of the basic examples of
singular lagrangian submanifolds, in particular those for which
results on their deformation spaces are available.
We start with the easiest case, a plane curve $C$ in $\dK^2$, given
as the zero set of a mapping $f:\dK^2 \rightarrow \dK$. Such a curve $C$
is obviously lagrangian. In this case the complex $\cC_C^\bullet$ is
simplifies to
$$
\begin{array}{rcl}
\cC_C^0 = \cO_C & \stackrel{\delta}{\longrightarrow} & \cC^1_C={\cH}\!om_{\cO_C}(\cI/\cI^2,\cO_C)
= {\cH}\!om_{\cO_C}(\cO_C,\cO_C) = \cO_C\\
h & \longmapsto & \{h,f\}
\end{array}
$$
It follows immediately that $\cH^2(\cC^\bullet_C) = 0$, while
$\cL\cT^1_C=\cH^1(\cC^\bullet_C)={\cC}\!oker(\delta)$. This sheaf
is supported on the singular points of the curve, let $x_0$ be such a point.
Then we have
$$
\cL\cT^1_{C,x_0}=\frac{\cO_{C,x_0}}{\{\{h,f\}\,|\,h\in\cO_{C,x_0}\}}
$$
Now the following equalities hold
\begin{eqnarray*}
\frac{\cO_{C,x_0}}{\left\{\{h,f\}| \, h \in \cO_{C,x_0}\right\}}
&=&
\frac{\Omega^2_{\dK^2,x_0}}{\left\{f\Omega^2_{\dK^2,x_0}+\{df \wedge dh|h \in \cO_{C,x_0}\}\right\}} \\
&=&
\frac{\Omega^2_{\dK^2,x_0}}{\left\{f\Omega^2_{\dK^2,x_0}+df\wedge d \Omega^0_{\dK^2,x_0}\right\}}
\end{eqnarray*}
because $\cO_{C,x_0} \cong \Omega^2_{\dK^2,x_0}/(f\Omega^2_{\dK^2,x_0})$ and the
Poisson bracket of two functions $f$ and $g$
corresponds under the isomorphism $\cO_{\dK^2,x_0} \cong \Omega_{\dK^2,x_0}$
to the $2$-form $df \wedge dg$. But it is known (see \cite{Malgrange}) that the dimension
of the last quotient equals $\mu$, the Milnor number of the plane curve
singularity $(C,x_0)$. So the result is:
$$
\cL\!\cT^1_C = \prod_{x_0 \in Sing(C)} \dK^{\mu(C,x_0)}
$$
This is remarkable because the usual $T^1_C$ has dimension $\tau$
(the Tjurina number) which is in general smaller than $\mu$. The
difference corresponds to the space of deformations of the restriction
of the symplectic structure to $L$ (see also \cite{Giv3}).

Applying lemma \ref{lemDecomp}, we see that the dimension of
$\cL\!\cT^1$ for a surface singularity which is a curve germ,
crossed with a smooth factor is also equal to the Milnor number of
this curve. This result can also be obtained by a direct calculus,
e.g., for a cuspidal edge given in four-space (with coordinates
$A,B,C,D$ and symplectic form $d\,A\wedge d\,C+d\,B\wedge d\,D$)
by the two equations $A, B^2-C^3$, we get $LT^1 = \dK^2$ and
$LT^2 = 0$.

We will proceed with further examples of lagrangian surfaces in $\dK^4$, which satisfy
condition P of theorem \ref{theFiniteness}. So
there are three strata: one point with embedding dimension four (supposed to
be the origin), the singular locus away from this point and the regular locus. In
order to simplify the calculation of the cohomology of
$\cC^\bullet$, we will suppose that our varieties
are strongly quasi-homogeneous in the sense of \cite{Mond}, that is,
one can choose local coordinates of the ambient space around each point of
$L$ such that the defining equations become weighted homogeneous
with positive weights. In this case, the \emph{de Rham}-complex
is a resolution of the constant sheaf as one can see by considering
the decomposition of the modules $\Omega^p_L$ into eigenspaces
of the Lie-derivative.
\begin{lemma}
Let $L\subset M$ be a strongly quasi-homogeneous lagrangian subvariety.
Consider the map $J:(\Omega^\bullet_L,d,\wedge)
\rightarrow (\cC^\bullet_L,\delta,\wedge)$ of DGA's from
corollary \ref{corJDach}. Denote by $\widetilde{\Omega}_L^\bullet$
the subcomplex ${\cI}\!m(J)$ in $\cC^\bullet_L$. Then
$\widetilde{\Omega}_L^\bullet$ is a resolution of $\dK_L$.
\end{lemma}
\begin{proof}
By the long exact cohomology sequence, it suffices to prove that
the complex ${\cK}er(J)$ is acyclic. This can be done in exactly
the same way as for $\Omega^\bullet_L$ provided that the inner
derivative $i_E$ ($E$ being the quasi-homogeneous Euler vector field)
maps ${\cK}er(J) \cap \Omega_L^p$ into ${\cK}er(J) \cap
\Omega_L^{p-1}$. But this follows from lemma \ref{lemTors}
because if $\omega$ is a torsion element than the same holds for
$i_E \omega$.
\end{proof}
\begin{korollar}
Denote by $\cG^\bullet_L$ the cokernel of the map $J$. Then there
is an exact sequence of $\cO_L$-modules
$$
0 \longrightarrow \widetilde{\Omega}_L^\bullet
  \longrightarrow \cC_L^\bullet
  \longrightarrow \cG_L^\bullet
  \longrightarrow 0
$$
and the long associated long exact sequence gives
$$
\cH^i(\cC_L^\bullet) = \cH^i(\cG_L^\bullet)
$$
for all $i \geq 0$. In particular, if $L$ is of dimension two, then we get
\begin{eqnarray*}
\cH^1(\cC^\bullet_L) &=& {\cK}\!er(\delta: \cG^1_L \rightarrow\cG^2_L) \\
\cH^2(\cC^\bullet_L) &=& {\cC}\!oker(\delta: \cG^1_L \rightarrow \cG^2_L)
\end{eqnarray*}
\end{korollar}
We can thus calculate $\cL\!\cT^1_L$ and $\cL\!\cT^2_L$ by
computing the induced morphism $\delta:\cG_L^1\rightarrow\cG^2_L$.
As $J$ is an isomorphism at smooth points, the sheaves $\cG_L^i$
are supported on the singular locus of $L$, which is of dimension
one. In a neighborhood of all of its regular points $q$
(points with embedding dimension three), the
germ is decomposable and the dimension of $\cH^i(\cC^\bullet_L)_q$
is given by lemma \ref{lemDecomp}. So we are only interested in
the one special point with maximal embedding dimension.
We now choose an element $p \in \cO_L$ which is finite when
restricted to the support of $\cG^i_L$, note that although this is
set-theoretically equal to the singular locus of $L$, it may have
embedded components. We will suppose that $p$ maps the origin in
$\dK^4$ to the origin in $\dK$. Consider the sheaves $p_*\cG_L^1$ and
$p_*\cG_L^2$, these are modules over $\cO_\dK$. Denote by
$\widetilde{E}$ resp. $\widetilde{F}$ the modules of section of
$p_*\cG_L^1$ resp. $p_*\cG_L^2$ in a small neighborhood of the
origin. Then they can be decomposed into torsion and torsion free
parts, the former being supported on the origin while the latter
is free over $\dK\{t\}$. In practice, this is done as follows:
As $\cG_L^1$ and $\cG^2_L$ are graded modules over $\cO_L$ and the
map $\delta:\cG_L^1\rightarrow\cG^2_L$ is homogeneous, we consider
the decomposition of these modules into homogeneous parts. The map $p$
is finite, so the torsion submodules of $\widetilde{E}$ and
$\widetilde{F}$ corresponds to homogeneous parts of $\cG_L^1$ and
$\cG^2_L$ in a \textbf{finite} number of degrees. This yields a
decomposition of $\widetilde{E}$ and $\widetilde{F}$ into
$\widetilde{E}=\widehat{E}\oplus E$ and $\widetilde{F}=
\widehat{F} \oplus F$ such that $\widehat{E}$ and $\widehat{F}$
are supported on the origin, while $E$ and $F$ are free.
$\widehat{E}$ and $\widehat{F}$ being artinian, the kernel and
cokernel of $\delta_{|\widehat{E}}$ can be computed explicitly.
The following lemma is used to do this for $\delta_{|E}$.
\begin{lemma}
The rank of $E$ and $F$ is the Milnor number $\mu$ of the
transversal curve singularity, i.e. the germ $(L',0)$ such that
$(L,p)$ = $(L',0) \times (\dK,0)$ for all $p \in \mbox{Sing}(L)
\,\backslash\, 0$. Therefore, $\delta_{|E}:E \rightarrow F$ is an
$(E,F)$-connection in the sense of \cite{Malgrange}.
\end{lemma}
\begin{proof}
This is an explicit calculation involving the definition of the
complex $\cC^\bullet_L$ and the map $J:\Omega^\bullet_L
\rightarrow \cC^\bullet_L$. It suffices to calculate the rank of
$(\cG^1_L)_p$ and $(\cG^2_L)_p$. So suppose that $(L,p)$ is a
decomposable germ. We choose coordinates $(x,y,s,t) \in \dK^4$ (with
symplectic form $\omega=dx \wedge dy + ds \wedge dt$) around $p$
such that $L$ is given as the zero locus of $s$ and a function
$f$ depending only on $x$ and $y$. Denote the ideal generated by
these two functions by $I$ and by $R$ the stalk of $\cO_L$ at the
point $p$. Then we can identify $I/I^2$ with $R^2$, so
$Hom_R(I/I^2,R)$ is free on the two generators $n_1$ and $n_2$,
where
$$
\begin{array}{cc}
n_1(f) = 1& n_1(s) = 0 \\
n_2(f) = 0& n_2(s) = 1
\end{array}
$$
while $Hom_R(I/I^2 \wedge I/I^2,R)$ is just $R$, generated by
the homomorphism sending $f\wedge s$ to $1$ in $R$.
The complex $\cC^\bullet$ at the point $p$ then reads:
$$
\begin{array}{ccccc}
R & \longrightarrow & R\,n_1 \oplus R\,n_2 & \longrightarrow & R \\
h & \longmapsto & (\{h,f\},\{h,s\})  \\
  &             & (p,q)                   & \longmapsto &  \{p,s\}+\{f,q\}
\end{array}
$$
where the pair $(p,q) \in R^2= Hom_R(I/I^2,R)$ denotes the
homomorphism sending $f \in I/I^2$ to $p \in R$ and $s \in I/I^2$
to $q \in R$.

Now we have to investigate the modules of differential forms on
$L$ at $x$. In general
$$
\Omega^p_R =  \Omega^p_S / (I \Omega^p_S + dI \wedge \Omega_S^{p-1})
$$
where $S$ is the ring $\dK\{x,y,s,t\}$. This leads to
\begin{eqnarray*}
\Omega^1_R & = & M_1 \oplus M_2 \\
\Omega^2_R & = & M_3 \oplus M_4
\end{eqnarray*}
where we have used the following abbreviations:
\begin{eqnarray*}
M_1 & = & \frac{R\,dx \oplus R\,dy}{R\,df} \\
M_2 & = & R\,dt \\
M_3 & = & \frac{R\,dx \wedge dy}{R\,df \wedge dx \oplus R\,df \wedge dy} \\
M_4 & = & \frac{R\,dx \wedge dt \oplus R\,dy \wedge dt }{R\,df \wedge dt}
\end{eqnarray*}
$J:\Omega^\bullet_L \rightarrow \cC^\bullet_L$ can be described as
\begin{eqnarray*}
J : M_1 & \longrightarrow & R\,n_1 \oplus R\,n_2 \\
dx & \longmapsto & \left(\{x,f\},\{x,s\}\right) = (\partial_y f,0)\\
dy & \longmapsto & \left(\{y,f\},\{y,s\}\right) = (-\partial_x f,0)\\
J : M_2 & \longrightarrow & R\,n_1 \oplus R\,n_2 \\
dt & \longmapsto & \left(\{t,f\},\{t,s\}\right) = (0,1)
\end{eqnarray*}
\begin{eqnarray*}
J : M_3 & \longrightarrow & R\\
dx\wedge dy& \longmapsto & J(dx)\wedge J(dy)=0 \\
J : M_4 & \longrightarrow & R\\
dx\wedge dt & \longmapsto & J(dx)\wedge J(dt) = \partial_y f\\
dy\wedge dt & \longmapsto & J(dx)\wedge J(dt) = -\partial_x f
\end{eqnarray*}
$E$ and $F$ are the cokernels of the maps $J:M_1\oplus M_2
\rightarrow R\,n_1\oplus R\,n_2$ and $J:M_3\oplus M_4
\rightarrow R$, respectively. So the result is
$$
E = F = R / \left(\partial_x f, \partial_y f \right)  =
\dK\{t\} \otimes \cO_{L',p}/ \left(\partial_x f, \partial_y f \right)  =
\dK\{t\}^\tau
$$
As $L$ is strongly quasi-homogeneous, we have weighted homogeneous local
equations for the transversal slice which gives $\tau = \mu$.
\end{proof}
Denote $\delta_{|E}$ by $D$ for short. Then $D$ is a first-order
differential operator $D:\cO_\dK^\mu \rightarrow \cO_\dK^\mu$
which respects the grading. So it is of the form
$$
D = t \partial_t \dI + A
$$
where $A$ is a constant $\mu \times \mu$-matrix. Thus, the second part
of the cohomology of $\cC_L^\bullet$ (i.e. kernel and cokernel of
$\delta$) can be deduced from the solutions of the differential
system given by $D$. All explicit calculations have been done using
\emph{Macaulay2}.

\begin{figure}
\centering \epsfig{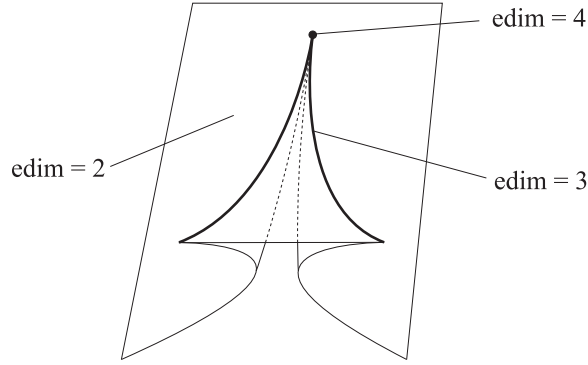}
\caption{The open swallowtail $\Sigma_2 \subset
\dK^4$} \label{sigma2}
\end{figure}

The first interesting example we are going to study is the so called ``open swallowtail''.
For details of its definition, see \cite{Giv1} and \cite{Giv2}.
Consider the space of polynomials in one variable of degree
$d:=2k+1$ with fixed leading coefficient and sum of roots equal to zero, that is, the space
$$
\begin{array}{rcccl}
\cP_{2k+1} & = &
\left\{x^{2k+1} + A_2 x^{2k-1}
+ \ldots + A_{2k+1}x^0 \right\} & \cong & \dK^{2k}\\
\end{array}
$$
which comes equipped with the following symplectic structure
$$
\omega =
\sum_{i=2}^{k+1} \left(2k+1-i\right)!
\left(i-2\right)! \cdot
(-1)^i dA_{i} \wedge dA_{2k+3-i}
$$
We will write $\Sigma_k$ for the subspace consisting of those
polynomials which have a root of multiplicity greater than
$k$. This space is obviously of dimension $k$ and it can be
shown that the form $\omega$ vanishes on its regular locus. So we
have a lagrangian subvariety in the space $\cP_{2k+1}$, which is
called open swallowtail. To get a more concrete impression of how it looks like,
we will describe the easiest examples. For $k=1$, $\Sigma_1 \subset \cP_{3}$ is just
the ordinary cusp in the plane, this case has already been discussed above. For
$k=2$, we obtain a surface in the four-dimensional space (see
the conceptual figure \vref{sigma2})
$$
\cP_5 =
\left\{
x^5+Ax^3+Bx^2+Cx+D \,|\, (A,B,C,D) \in \dK^4
\right\}
$$
(the symplectic form is $\omega=3d\,A\wedge d\,D+d\,C\wedge d\,B$)
consisting of those polynomials $f$ with a root of multiplicity at least three.
Such a $f$ can be written as $f= (x-a)^3(x^2+3ax+b)$, so there is a normalization of $\Sigma_2$ given by
\begin{eqnarray*}
n:\dK^2 &\longrightarrow& \cP_5 = \dK^4\\
(a,b) & \longmapsto & (b-6a^2,8a^3-3ab,3a^2b-3a^4,-a^3b)
\end{eqnarray*}
Note that the singular locus of $\Sigma_2$ is a again a cusp as
well as the transversal curve singularity.

The space $\Sigma_2$ is our main example, we will describe in some
detail how to apply the general results in this case.
Using elimination theory, we can calculate the defining equations of $\Sigma_2$
in $\dK^4$. It turns out that the swallowtail is a determinantal variety given
by the minors of the matrix
$$
\begin{pmatrix}
9D& 9 B^2-32 AC\\
3 C& -5 AB+ 125 D\\
-9 B& 45 A^2-100C
\end{pmatrix}
$$
The ideal which defines $\Sigma_2$ is generated by the following
three polynomials
\begin{eqnarray*}
f_1 & = &  -27 B^2 C+96 A C^2-45 A B D+1125 D^2, \\
f_2 & = &  81 B^3-288 A B C+405 A^2 D-900 C D \\
f_3 & = & -45 A B^2+135 A^2 C-300 C^2+1125 B D
\end{eqnarray*}
So $\Sigma_2$ is not a complete intersection but nevertheless
Cohen-Macaulay by the Hilbert-Burch theorem. We list
the commutators $\{f_i,f_j\}$ (for $1\leq i<j \leq 3$)
with respect to the given set of generators (this is a direct proof that
$\Sigma_2 \subset \dK^4$ is involutive):
\begin{eqnarray*}
\{f_1,f_2\} & = & -576A f_1 + 81B f_2 - 96 C f_3  \\
\{f_1,f_3\} & = &  15 A f_2 -12 B f_3  \\
\{f_2,f_3\} & = & -900 f_1 + 18 A f_3   \\
\end{eqnarray*}
$\Sigma_2$ is quasi-homogeneous with the
weights $(2,3,4,5)$ for the variables $A$,
$B$, $C$, $D$, respectively. We can thus apply the machinery
developed above to obtain that $\left(\cL\!\cT^1_{\Sigma_2}\right)_0=0$,
while $\left(\cL\!\cT^2_{\Sigma_2}\right)_0=\dK$.
The operator $D$ is in this case
$$
t\partial_t\,\dI+
\begin{pmatrix}
11/40& -245/2& 0& 0\\
33/4000& 109/40& 0& 0\\
0& 0& 49/15& -59/27\\
0& 0& 51/100& 11/15
\end{pmatrix}
$$
For $\dK=\dC$, the monodromy of the locally constant
sheaf $\cL\!\cT^1_{|Sing(\Sigma_2)\,\backslash 0\,}$
has the following eigenvalues
$$
-\frac{8}{10}
\,,\,
-\frac{13}{10}
\,,\,
-\frac{22}{10}
\,,\,
-\frac{27}{10}
$$

The second large class of examples are the conormal spaces. Given any
submanifold $Y$ of an $n$-dimensional manifold $X$, the total space of the
conormal bundle $T_Y^* X$ is always a lagrangian submanifold of $T^*X$. More
generally, if $Y$ is an analytic subspace, we can take the closure of the
space of conormals to all smooth points of $Y$. The result (which is called
conormal space of $Y$ in $X$) is still lagrangian, but may have singularities.
This is an important class of lagrangian subvarieties,
as the characteristic variety of a holonomic $\cD_X$-module is always a finite union
of conormal spaces. Obviously, these spaces are conical in the fibers of
$T^*X$. If $X$ is a plane curve in $C \subset \dK^2$, then the conormal
space $T^*_C \dK^2$ will be a surface in $\dK^4$. Here the results
are as follows.

\renewcommand{\baselinestretch}{1.4} \small
$$
\begin{array}{r|l|l|l}
\textbf{\mbox{equation of $C$}}& \mathbf{LT^1} &
\mathbf{LT^2} & \textbf{\mbox{eigenvalues (multiplicity, if $\neq 1$)}} \\
\hline
y^2-x^5 & 0 & 0 &  -\frac{4}{5},-\frac{16}{5}\\
\hline
y^3-x^7 & 0 & 0 &  -\frac{37}{7},-\frac{61}{7},-\frac{69}{7},-\frac{85}{7},
-\frac{93}{7},-\frac{117}{7}\\
\hline
y^5-x^7 & 0 & 0 &  -\frac{116}{7},-\frac{132}{7},-\frac{148}{7},-\frac{164}{7},\\
\hline
y^3-x^6 & \dK & \dK & -\frac{7}{2},-\frac{10}{2}^{(2)},-\frac{13}{2}\\
\hline
x y (x+y)(x-y)(x-2y)& \dK^2 & \dK^2 & - \\
\end{array}
$$
\renewcommand{\baselinestretch}{1}
\normalsize

\nd
In the last example, there is only an isolated singularity, so the
modules $\cG_L^1$ and $\cG_L^2$ are artinien.

Finally, there is a third class of singular lagrangian
subvarieties, these are \emph{completely integrable hamiltonian
systems}. Such a system is given in the $2n$-dimensional phase
space by $n$ Poisson-commuting functions. The ideal formed by them
then obviously satisfies the involutivity condition. If,
additionally, the common zero set of these function is a complete
intersection, then it will be lagrangian in our sense.
The lagrangian deformation space of such a system is at least
$n$-dimensional (addition of a constant is flat and the ideal
stays involutive).

To get the equations of some interesting examples, we will proceed as
follows. Choose coordinates $(p_1,q_1,p_2,q_2)$ of $\dK^4$ and set
$z_1 = p_1 + i q_1$ and $z_2 = p_2 + i q_2$ (This can obviously be done
only in the real case, but it is a formal calculus which works as well for
$\dK=\dC$ as for $\dK=\real$). We can now express functions on $\dK^4$ in the
variables $z_1, z_2, \overline{z_1}, \overline{z_2}$, and the Poisson bracket
becomes
$$
\{f,g\} =
2i \left(
  \partial_{\overline{z}_1}f \cdot \partial_{z_1}g
- \partial_{\overline{z}_1}g \cdot \partial_{z_1}f
+ \partial_{\overline{z}_2}f \cdot \partial_{z_2}g
- \partial_{\overline{z}_2}g \cdot \partial_{z_2}f
\right)
$$
We want to find functions $f_1, f_2$ such that $\{f_1,f_2\} = 0$. Set, for
example $f  =  \lambda z_1 \overline{z_1} + \mu z_2 \overline{z_2}$ and let us
look for a $g=  z_1^\alpha\overline{z_1}^\beta z_2^\gamma
\overline{z_2}^\delta$ for some parameters $\lambda,\mu,\alpha,\beta,\gamma,\delta \in\dN$.
It can be easily verified that the commuting condition transforms to
$$
\lambda(\alpha-\beta)-\mu(\gamma-\delta)=0
$$
The following table shows results for some resonance (\textbf{r}) coefficients
$\lambda,\mu$ and exponents  (\textbf{e})
$\alpha,\beta,\gamma,\delta$.

\renewcommand{\baselinestretch}{1.4} \small
$$
\begin{array}{r|c|l|l|l}
\textbf{\mbox{r}}& \textbf{\mbox{e}}&
\mathbf{LT^1} & \mathbf{LT^2} &
\textbf{\mbox{eigenvalues (multiplicity)}} \\
\hline
1 , 0 & 0,0,1,1 & \dK^2 & \dK &  -3^{(4)} \\
\hline
1 , 2 & 0,2,1,0 & \dK^3 & \dK^2 &-\frac{2}{2}^{(2)},-\frac{3}{2}^{(2)},-\frac{4}{2}^{(2)},
-\frac{5}{2}^{(2)},-\frac{6}{2}^{(2)} \\
\hline
1 , 3 & 3,0,0,1 & \dK^4 & \dK^3 &-\frac{3}{3}^{(2)},-\frac{5}{3}^{(2)},-\frac{7}{3}^{(4)},
-\frac{9}{3}^{(4)},-\frac{11}{3}^{(4)},-\frac{13}{3}^{(2)},-\frac{15}{3}^{(2)}\\
\hline
1 , 4 & 4,0,0,1 & \dK^5 & \dK^4 &  -\frac{4}{4}^{(2)},-\frac{7}{4}^{(2)},-\frac{9}{4}^{(2)},
-\frac{10}{4}^{(2)},-\frac{12}{4}^{(2)},-\frac{13}{4}^{(2)},\\
%\cline{1-4}
&&&& -\frac{14}{4}^{(2)},-\frac{15}{4}^{(2)},-\frac{16}{4}^{(2)},-\frac{17}{4}^{(2)},
-\frac{18}{4}^{(2)},-\frac{19}{4}^{(2)},\\
%\cline{1-4}
&&&&-\frac{20}{4}^{(2)},-\frac{22}{4}^{(2)},-\frac{23}{4}^{(2)},
-\frac{25}{4}^{(2)},-\frac{28}{4}^{(2)}\\
\end{array}
$$
\renewcommand{\baselinestretch}{1} \normalsize

\paragraph{\textbf{Remark:}} The eigenvalues in all examples
have a symmetry property, which we cannot prove at this moment.
These eigenvalues looks very similar to the spectrum of an
isolated hypersurface singularity. One might speculate that
there is a mixed Hodge structure related to this theory and
that the eigenvalues share further properties with the
spectrum, e.g. the semi-continuity under deformations.

%\bibliographystyle{amsalpha}
%\bibliography{Lagrange}

\providecommand{\bysame}{\leavevmode\hbox to3em{\hrulefill}\thinspace}

\vspace*{1cm}

\nd
FB 17, Mathematik,\\
Johannes-Gutenberg-Universit\"at Mainz,\\
D-55099 Mainz, Germany \\

\nd
straten@mathematik.uni-mainz.de,\\
sevenhec@mathematik.uni-mainz.de

\end{document}